\documentclass{article}
\usepackage{amssymb}
\usepackage{amsfonts}
\usepackage{amsmath}
\usepackage{geometry}
\usepackage[doublespacing]{setspace}

\newtheorem{theorem}{Theorem}

\newtheorem{corollary}[theorem]{Corollary}

\newtheorem{definition}[theorem]{Definition}
\newtheorem{example}[theorem]{Example}

\newtheorem{lemma}[theorem]{Lemma}

\newtheorem{remark}[theorem]{Remark}

\newenvironment{proof}[1][Proof]{\noindent\textbf{#1.} }{\ \rule{0.5em}{0.5em}}
\input{tcilatex}

\begin{document}

\title{Normal Quadratics in Ore Extensions, Quantum Planes, and Quantized
Weyl Algebras}
\author{Candis Holtz \\
N 5894 County BB\\
Cecil, WI 54111 \and Kenneth Price \\
University of Wisconsin Oshkosh\\
Oshkosh, WI 54901}
\date{June 6, 2008}
\maketitle

\begin{abstract}
Let $R$ be a Noetherian domain and let $\left( \sigma ,\delta \right) $ be a
quasi-derivation of $R$ such that $\sigma $ is an automorphism. \ There is
an induced quasi-derivation on the classical quotient ring $Q$ of $R$. \
Suppose $F=t^{2}-v$ is normal in the Ore extension $R\left[ t;\sigma ,\delta %
\right] $ where $v\in R$. \ We show $F$ is prime in $R\left[ t;\sigma
,\delta \right] $\ if and only if $F$ is irreducible in $Q\left[ t;\sigma
,\delta \right] $ if and only if there does not exist $w\in Q$ such that $%
v=\sigma \left( w\right) w-\delta \left( w\right) $. \ We apply this result
to classify prime quadratic forms in quantum planes and quantized Weyl
algebras.
\end{abstract}

\section{Overview\label{Overview Sec}}

In a ring $R$ a nonzero element $a$ is normal if $aR=Ra$. \ Moreover, if $R$
is prime then $a$ induces an automorphism $\varphi $ such that $ar=\varphi
\left( r\right) a$ for all $r\in R$. We would like to know as much as
possible about\ normal elements in Ore extensions. \ Jordan studied degree
one normal elements of Ore extensions and applied his results to quantum
determinants in $n\times n$ quantum matrices. \ In section \ref{Normal Ore
Sec}, we prove Lemma \ref{Normal Lemma}, which was motivated by his proof of 
\cite[Proposition 1]{Jordan}. \ This is used in Theorem \ref{Factor Theorem}
to classify prime monic normal polynomials of degree 2 in an Ore extension.

Theorem \ref{Factor Theorem} is used to help establish the main result of 
\cite{Price}, which gives a test to determine when the universal enveloping
algebra of a Lie color algebra is a domain. \ To find new applications of
Theorem \ref{Factor Theorem} we study factorization in quantum planes and
quantized Weyl algebras. \ Background is provided in section \ref{QP & QWA
sec}.

One of the coauthors, Holtz, first encountered the topic of factorization in
rings while enrolled in an undergraduate algebra course taught by the
faculty coauthor, Price. \ We pursued ways to extend results from
commutative polynomial rings to quantized Weyl algebras. \ Similar results
on quantum planes (which are extended here) were completed by the faculty
coauthor and another undergraduate student (see \cite{CP}). \ Mrs. Holtz
assisted with the results in section \ref{QWA Def sec}.

We apply Theorem \ref{Factor Theorem} to quantum planes and quantized Weyl
algebras in section \ref{prime sec}. \ An easy degree argument proves prime
polynomials are irreducible. \ In sections \ref{QP & QWA sec} and \ref{QWA
Def sec} we derive reducibility results which are needed to determine the
prime quadratic forms in a quantum plane (Theorem \ref{Quad Forms in QP})
and in a quantized Weyl algebra (Theorem \ref{Main Theorem}).

Throughout $k$ denotes a commutative field of characteristic different from
2 and\ all rings are equipped with a multiplicative identity. \ We refer to 
\cite[Section 1.4]{Bueso} for background on factorization in noncommutative
domains. \ Our approach to factorization is two-sided, so we cover basic
concepts.

\begin{definition}
Suppose $R$ is a domain and $a\in R$.

\begin{enumerate}
\item We call $a$ \emph{irreducible} if $a$ is not a unit and if $a=bc$ for
some $b,c\in R$ then either $b$ is a unit or $c$ is a unit.

\item We say $a$ \emph{divides} $c\in R$, and write $a|c$, if there exists $%
b\in R$ such that $c=ab$ or $c=ba$. \ In this case we call $c$ a \emph{%
multiple} of $a$.

\item We call $a$ \emph{prime} if $a|bc$ implies $a|b$ or $a|c$ for all $%
b,c\in R$.
\end{enumerate}
\end{definition}

It is easy to see a normal element of $R$ is prime if and only if the
principal ideal it generates is completely prime. \ Thus classifying normal
elements as prime also establishes existence of completely prime ideals.

\section{Monic Normal Polynomials of Degree 2 in $R\left[ t;\protect\sigma ,%
\protect\delta \right] $ \label{Normal Ore Sec}}

Throughout this section $R$, is a Noetherian domain and $\sigma $ is an
automorphism of $R$.\ \ Recall that a $\sigma $-derivation is an additive
map $\delta :R\rightarrow R$ such that $\delta \left( rs\right) =\delta
\left( r\right) s+\sigma \left( r\right) \delta \left( s\right) $ for all $%
r,s\in R$. \ The pair $\left( \sigma ,\delta \right) $ is called a
quasi-derivation of $R$. \ Set $S=R\left[ t;\sigma ,\delta \right] $, the
Ore extension over $R$, which is a Noetherian domain by \cite[Theorem 2.6,
p. 39; and Exercise 2.2O, p. 37]{GW}.

\begin{lemma}
\label{Normal Lemma}Suppose $F\,=t^{2}-v$ is normal in $R\left[ t;\sigma
,\delta \right] $ where $v\in R$.\ Then $v$ commutes with $t$ and the
following hold for all $r$ in $R$.

\begin{enumerate}
\item $\delta ^{2}\left( r\right) =vr-\sigma ^{2}\left( r\right) v$

\item $\left( \delta \sigma \right) \left( r\right) =-\left( \sigma \delta
\right) \left( r\right) $

\item $F\,r=\sigma ^{2\,}\left( r\right) \,F$ for all $r\in R$
\end{enumerate}
\end{lemma}

\begin{proof}
Set $S=R\left[ t;\sigma ,\delta \right] $ and let $\varphi $ be the
automorphism of $S$ such that $F\,s=\varphi \left( s\right) \,F\,$\ for all $%
s\in S$. \ Then $\varphi $ preserves degree in $t$ since $S$ is a domain. \
Therefore $\varphi \left( t\right) =at+b$ for some $a,b\in R$ with $a\neq 0$%
. \ We obtain equation \ref{first eq}\ by expanding both sides of $%
F\,t=\varphi \left( t\right) \,F$. \ We find $v$ commutes with $t$ by
equating coefficients in equation \ref{first eq}.%
\begin{equation}
t^{3}-vt=at^{3}+bt^{2}-a\sigma \left( v\right) t-\left( a\delta \left(
v\right) +bv\right)  \label{first eq}
\end{equation}

Let $r\in R$ be arbitrarily chosen. \ We expand both sides of $\varphi
\left( r\right) F=Fr$ to obtain equation \ref{commutator equation}. \ Expand
and equate coefficients in equation \ref{commutator equation} to obtain 1,
2, and 3.%
\begin{equation}
\varphi \left( r\right) t^{2}-\varphi \left( r\right) v=\sigma ^{2}\left(
r\right) t^{2}+\left( \delta \sigma +\sigma \delta \right) \left( r\right)
t+\delta ^{2}\left( r\right) -vr  \label{commutator equation}
\end{equation}
\end{proof}

We may form the classical quotient ring $Q$ of $R$, which is a division
ring. \ Following \cite[Lemma 1.3]{Good}, we may extend $\left( \sigma
,\delta \right) $ to a quasi-derivation of $Q$. \ Then $Q\left[ t;\sigma
,\delta \right] $ is a principal left and right ideal domain by \cite[%
Theorem 2.8, p. 39]{GW}.

\begin{theorem}
\label{Factor Theorem}Set $F\,=t^{2}-v$ with $v\in R$ and let $S=R\left[
t;\sigma ,\delta \right] $.

\begin{enumerate}
\item If $F$ is reducible in $Q\left[ t;\sigma ,\delta \right] $ then $%
F=\left( t-\sigma \left( w\right) \right) \left( t+w\right) $ for some $w\in
Q$. \ In this case $v=\sigma \left( w\right) w-\delta \left( w\right) $.

\item If $F$ is normal in $S$ then $F$ is prime in $S$ if and only if $F$ is
irreducible in $Q\left[ t;\sigma ,\delta \right] $.
\end{enumerate}
\end{theorem}

\begin{proof}[Proof]
If $F$ is reducible then there is a factorization of the form in equation %
\ref{factorization} for some $a_{0},a_{1},b_{0},b_{1}\in Q$ with $a_{0}$ and 
$b_{0}$ both nonzero. 
\begin{equation}
F=\left( a_{0}t+a_{1}\right) \left( b_{0}t+b_{1}\right)
\label{factorization}
\end{equation}%
\ Set $w=\left( b_{0}\right) ^{-1}b_{1}$. \ We can rewrite equation \ref%
{factorization} as equation \ref{rewrite}.%
\begin{equation}
t^{2}-v=\left( a_{0}\sigma \left( b_{0}\right) t+a_{0}\delta \left(
b_{0}\right) +a_{1}b_{0}\right) \left( t+w\right)  \label{rewrite}
\end{equation}%
By expanding the right hand side of equation \ref{rewrite} and equating
coefficients we obtain $a_{0}\sigma \left( b_{0}\right) =1$ and $\sigma
\left( w\right) =-\left( a_{0}\delta \left( b_{0}\right) +a_{1}b_{0}\right) $%
. \ Part 1 follows immediately.

Let $I$ be the principal ideal of $S$ generated by $F$. \ It is easy to see $%
Q\left[ t;\sigma ,\delta \right] I$ is completely prime if and only if $F$
is irreducible in $Q\left[ t;\sigma ,\delta \right] $ since $Q$ is a
division algebra. \ If $Q\left[ x;\sigma ,\delta \right] I\cap S=I$ then the
equivalence follows by \cite[Lemma 8.1.21]{Bueso}. \ Let $P_{1}\in Q\left[
x;\sigma ,\delta \right] I\cap S$ be arbitrarily chosen. \ By \cite[%
Proposition 2.1.16 (iv)]{McRob} there exists $a\in R\backslash \left\{
0\right\} $ such that $aP_{1}\in I$. \ Therefore $aP_{1}=P_{2}F$ for some $%
P_{2}\in S$. \ We set $m_{i}=\deg P_{i}$ and write $P_{i}$ as in equation %
\ref{psum} with $r_{i,j}\in R$ for $j\in \left\{ 0,\ldots ,m_{i}\right\} $
and $i=1,2$. 
\begin{equation}
P_{i}=\sum_{j=0}^{m_{i}}r_{i,j}t^{j}  \label{psum}
\end{equation}%
We have $vt^{j}=t^{j}v$ for all $j$ by Lemma \ref{Normal Lemma}. \ Thus by
expanding both sides of the equation $aP_{1}=P_{2}F$ and equating
coefficients we obtain $m_{1}=m_{2}+2$, $r_{2,m_{2}}=ar_{1,m_{2}+2}$, $%
r_{2,m_{2}-1}=ar_{1,m_{2}+1}$, and the remaining coefficients satisfy the
recursive formula in equation\ \ref{recursive form}\ for $j=2,\ldots ,m_{2}$%
. \ 
\begin{equation}
r_{2,j-2}=ar_{1,j}+r_{2,j}v  \label{recursive form}
\end{equation}

Set $r_{3,m_{2}}=r_{1,m_{2}+2}$ and $r_{3,m_{2}-1}=r_{1,m_{2}+1}$. \ A
straightforward induction argument using equation \ref{recursive form}
proves there exists $r_{3,j}\in R$ such that $r_{2,j}=ar_{3,j}$ for $%
j=0,\ldots ,m_{2}$. \ Set $m_{3}=m_{2}$ and define $P_{3}$ as in equation %
\ref{psum} with $i=3$. \ Then $P_{1}=P_{3}F$ since $P_{2}=aP_{3}$ by
construction. \ Therefore $P_{1}\in I$ and, since $P_{1}$ was arbitrarily
chosen, $Q\left[ x;\sigma ,\delta \right] I\cap S=I$ as desired.
\end{proof}

\section{Quantum Planes and Quantized Weyl Algebras \label{QP & QWA sec}}

Let $q\in k\backslash \left\{ 0\right\} $ be arbitrary. \ There is an
automorphism of $k\left[ x\right] $ given by $\sigma \left( x\right) =qx$. \
If $q\neq 1$ then by \cite[Exercise 2M, p. 33]{GW} there is a $\sigma $%
-derivation $\delta $ of $k\left[ x\right] $ such that 
\begin{equation}
\delta \left( P\right) =\dfrac{\sigma \left( P\right) -P}{\left( q-1\right) x%
}  \label{eulerian derivative}
\end{equation}%
for all $P\in k\left[ x\right] $. \ In particular, $\delta \left( x\right)
=1 $.

The quantum plane is the Ore extension denoted $\mathcal{O}_{q}\left(
k^{2}\right) =k\left[ x\right] \left[ y;\sigma \right] $ and the quantized
Weyl algebra is the Ore extension $\mathcal{A}_{1}^{q}\left( k\right) =k%
\left[ x\right] \left[ y;\sigma ,\delta \right] $. \ In case $q=1$ we have $%
\mathcal{O}_{1}\left( k^{2}\right) =k\left[ x,y\right] $, the ordinary
commutative polynomial ring in two variables, and $\mathcal{A}_{1}^{1}\left(
k\right) =\mathcal{A}_{1}\left( k\right) $, the ordinary Weyl algebra. \ See 
\cite{Co} for details on how to obtain the Weyl algebra as a ring of
differential operators. \ If $q$ is not a root of unity then $\mathcal{A}%
_{1}^{q}\left( k\right) $ can be constructed in a similar way if the
ordinary derivative is replaced with the one in equation \ref{eulerian
derivative}.

\begin{remark}
\label{Center Desc}Suppose $q$ is a primitive $n^{\text{th}}$ root of unity
and $P\in \mathcal{A}_{1}^{q}\left( k\right) $ or $P\in \mathcal{O}%
_{q}\left( k^{2}\right) $. Then $P$ is central if and only if the exponent
of $x$ and the exponent of $y$ are both multiples of $n$ for every term in
the support of $P$. \ This is well-known (see for example \cite[Lemma 2.2]%
{AVDVO}).
\end{remark}

\begin{definition}
A \emph{quadratic form} is an element $F$ of $\mathcal{O}_{q}\left(
k^{2}\right) $ or $\mathcal{A}_{1}^{q}\left( k\right) $ which can be written
as in equation \ref{QF Eq} for some $a,b,c,d,e,f\in k$ with at least one of $%
a,b,c$ nonzero.%
\begin{equation}
F\left( x,y\right) =ax^{2}+bxy+cy^{2}+dx+ey+f  \label{QF Eq}
\end{equation}%
The \emph{quantum discriminant} of $F$ is $\Delta _{q}\left( F\right)
=b^{2}-4acq$.
\end{definition}

If $F$ is reducible then it can be expressed as in equation \ref{reduce eq}
for some $\alpha ,\beta ,\gamma ,\lambda ,\mu ,\nu \in k$.%
\begin{equation}
F=\left( \lambda x+\mu y+\nu \right) \left( \alpha x+\beta y+\gamma \right)
\label{reduce eq}
\end{equation}%
Expanding and equating coefficients of degree two terms leads to $\Delta
_{q}\left( F\right) =\left( \lambda \beta -q\mu \alpha \right) ^{2}$. \ Thus
if $\Delta _{q}F$ is not a square then $F$ is irreducible in $\mathcal{O}%
_{q}\left( k^{2}\right) $ or $\mathcal{A}_{1}^{q}\left( k\right) $. \
Homogeneous irreducible quadratic forms in $\mathcal{O}_{q}\left(
k^{2}\right) $ were classified in \cite{CP}. \ Theorem \ref{Quad Forms in QP}
extends this slightly.

\begin{theorem}
\label{Quad Forms in QP}Let $F=ax^{2}+bxy+cy^{2}+f$ be a quadratic form in $%
\mathcal{O}_{q}\left( k^{2}\right) $. \ Then $F$ is reducible if and only if
one of the following holds.

\begin{enumerate}
\item $f=0$ and $\Delta _{q}F$ is a square in $k$

\item $f\neq 0$, there exists $\delta \in k$ such that $\Delta _{q}F=\delta
^{2}$ and $\delta \left( 1+q\right) =b\left( q-1\right) $, and either (a) $%
a\neq 0$ and $-af$ is a square in $k$ or (b) $c\neq 0$ and $-cf$ is a square
in $k$.
\end{enumerate}
\end{theorem}

\begin{proof}
Case 1 is handled by \cite[Theorem 6]{CP}. \ Case 2(b) may be handled the
same way as Case 2(a). \ We assume $a\neq 0$ and $f\neq 0$. \ If case 2(a)
holds then there exists $\gamma \in k$ such that $\left( a\gamma \right)
^{2}=-af$. \ Moreover $F$ is reducible with a factorization given in
equation \ref{case 2a}. \ 
\begin{equation}
F=\left( ax+\frac{b+\delta }{2q}y+\frac{f}{\gamma }\right) \left( x+\frac{%
b-\delta }{2a}y+\gamma \right)  \label{case 2a}
\end{equation}%
On the other hand, if $F$ is reducible then we may write it as in equation %
\ref{qp reduce}.%
\begin{equation}
F=\left( ax+\frac{b-a\beta }{q}y+\frac{f}{\gamma }\right) \left( x+\beta
y+\gamma \right)  \label{qp reduce}
\end{equation}%
Expanding and equating coefficients leads to $\left( a\gamma \right)
^{2}=-af $, $a\beta ^{2}-b\beta +qc=0$, and $b=a\beta \left( 1+q\right) $. \
The first equation is 2(a) and since the second equation has a solution we
know $\Delta _{q}\left( f\right) =\delta ^{2}$ for some $\delta \in k$ such
that $\beta =\left( b-\delta \right) /2a$. \ This gives $b\left( q-1\right)
=\delta \left( 1+q\right) $, as desired. \ 
\end{proof}

\section{Reducibility in $\mathcal{A}_{1}^{q}\left( k\right) $ \label{QWA
Def sec}}

Checking quadratic forms for reducibility in $\mathcal{A}_{1}^{q}\left(
k\right) $ is more complicated since the product of two polynomials may be
homogeneous even if the polynomials themselves are not homogeneous. \ For
example, $\left( x+y+1\right) \left( x+y-1\right) =x^{2}+\left( 1+q\right)
xy+y^{2}$.

\begin{theorem}
\label{Quad Form Thm}Let $F=ax^{2}+bxy+cy^{2}+f$ be a quadratic form in $%
\mathcal{A}_{1}^{q}\left( k\right) $. \ Then $F$ is reducible if and only if
one of the following holds.

\begin{enumerate}
\item either (a) $a=f=0$, or (b) $c=f=0$, or (c) $a=b=0$ and $-cf$ is a
square in $k$, or (d) $b=c=0$ and $-af$ is a square in $k$

\item $f\neq 0$ and $\Delta _{q}F=\left( b-2qf\right) ^{2}$

\item $q\neq -1$, $b\neq 0$, $\Delta _{q}F=\left( \frac{1-q}{1+q}b\right)
^{2}$, and $\sigma ^{2}=a\left( \frac{b}{1+q}-f\right) $ for some $\sigma
\in k$

\item $q=-1$, $b=0$, $\Delta _{-1}F=4\tau ^{2}\neq 0$, and $\left( \tau
-f\right) c=\omega ^{2}$ for some $\tau ,\omega \in k$
\end{enumerate}
\end{theorem}

\begin{proof}
First we assume $F$ is reducible and write $F$ as in equation \ref{reduce eq}%
. \ Expanding and equating coefficients yields six equations.%
\begin{equation*}
\begin{tabular}{rllrl}
\textrm{i} & $\lambda \alpha =a$ & \hspace{0.75in} & \textrm{iv} & $\lambda
\gamma +\nu \alpha =0$ \\ 
\textrm{ii} & $\lambda \beta +q\mu \alpha =b$ &  & \textrm{v} & $\mu \gamma
+\nu \beta =0$ \\ 
\textrm{iii} & $\mu \beta =c$ &  & \textrm{vi} & $\nu \gamma +\mu \alpha =f$%
\end{tabular}%
\end{equation*}

Suppose $\mu \alpha =0$. \ Then equations i and iii imply $a=0$ or $c=0$. \
If $b\neq 0$ and $f\neq 0$ then\ Equations ii and vi imply $\beta ,\gamma
,\lambda ,\nu $ are all nonzero. \ Then $\alpha \neq 0$ and $\mu \neq 0$ by
equations iv and v, which contradicts $\mu \alpha =0$. \ Thus we are left
with $b=0$ or $f=0$. \ If $f=0$ then case 1(a) or 1(b) holds. \ On the other
hand if $b=0$ and $f\neq 0$ then $F$ is a quadratic polynomial in
one-variable thus case 1(c) or 1(d) must hold.

Now suppose $\mu \alpha \neq 0$. \ Set $\pi =\mu \alpha $, $\chi =\nu \alpha 
$, and $\rho =\mu \gamma $. \ Then equation \ref{reduce eq} may be rewritten
as equation \ref{simpler eq}.%
\begin{equation}
F=\frac{1}{\pi }\left( ax+\pi y+\chi \right) \left( \pi x+cy+\rho \right)
\label{simpler eq}
\end{equation}%
Expanding and equating coefficients yields four equations, listed as vii-x
below.%
\begin{equation*}
\begin{tabular}{rllrl}
\textrm{vii} & $q\pi ^{2}+ac=b\pi $ & \hspace{0.75in} & \textrm{ix} & $\pi
\chi +a\rho =0$ \\ 
\textrm{viii} & $\chi \rho +\pi ^{2}=f\pi $ &  & \textrm{x} & $c\chi +\pi
\rho =0$%
\end{tabular}%
\end{equation*}

Equation ix gives $\chi =-\frac{a\rho }{\pi }$ which we substitute into viii
and x to find equations \ref{rho square eq} and \ref{rho eq}.%
\begin{equation}
f=\pi -a\left( \frac{\rho }{\pi }\right) ^{2}  \label{rho square eq}
\end{equation}%
\begin{equation}
\left( ac-\pi ^{2}\right) \rho =0  \label{rho eq}
\end{equation}%
We are left with three subcases.

\begin{itemize}
\item Suppose $\rho =0$. \ Then equations \ref{rho square eq} and vii imply $%
f=\pi $ and $b=qf+\frac{1}{f}ac$.\ \ The quantum discriminant is $\Delta
_{q}F=\left( b-2qf\right) ^{2}$ so case 2 holds.

\item Suppose $\rho \neq 0$ and $q\neq -1$. \ Equations \ref{rho eq} and vii
give $ac=\pi ^{2}$ and $b=\left( 1+q\right) \pi $ so $b\neq 0$ and $\Delta
_{q}F=\left( \frac{1-q}{1+q}b\right) ^{2}$. \ Substituting $b=\left(
1+q\right) \pi $ into equation \ref{rho square eq} gives $\left( \frac{a\rho 
}{\pi }\right) ^{2}=a\left( \frac{b}{1+q}-f\right) $. \ Thus case 3 holds.

\item Suppose $\rho \neq 0$ and $q=-1$. \ Equations \ref{rho eq} and vii
give $ac=\pi ^{2}$ and $b=0$ so $\Delta _{q}F=4\pi ^{2}\neq 0$. \ Equation %
\ref{rho square eq} leads to $\rho ^{2}=\left( \pi -f\right) c$. \ Thus case
4 holds.
\end{itemize}

Equations \ref{factorization 1}, \ref{factorization 3}, and \ref%
{factorization 2} provide factorizations of $F$ in case 2, 3, and 4. \ Case
1 is trivial%
\begin{equation}
F=\left( ax+fy\right) \left( x+\frac{1}{f}cy\right)  \label{factorization 1}
\end{equation}%
\begin{equation}
F=\left( ax+\tau y-\frac{a\omega }{\tau }\right) \left( x+\frac{1}{\tau }cy+%
\frac{\omega }{\tau }\right)  \label{factorization 3}
\end{equation}%
\begin{equation}
F=\left( ax+\frac{b}{1+q}y+\sigma \right) \left( x+\frac{1+q}{b}cy-\frac{1}{a%
}\sigma \right)  \label{factorization 2}
\end{equation}%
In cases 2 and 4 factorizations are given by equations \ref{factorization 1}
and \ref{factorization 3}, respectively. \ In case 3, the conditions imply $%
a\neq 0$. \ A factorization is given in equation \ref{factorization 2}.
\end{proof}

\begin{corollary}
\label{Quad Form Coro}Let $a,b,c\in k\backslash \left\{ 0\right\} $ and $%
f\in k$. \ The following hold in $\mathcal{A}_{1}^{q}\left( k\right) $.

\begin{enumerate}
\item $bxy+f$ is reducible if and only if $f=0$ or $b=qf$.

\item $ax^{2}+cy^{2}+f$ is reducible if and only if one of the following
holds

\begin{enumerate}
\item $f\neq 0$ and $ac=-qf^{2}$

\item $q=-1$, $ac=\tau ^{2}$ and $\left( \tau -f\right) c=\omega ^{2}$ for
some $\tau ,\omega \in k$
\end{enumerate}
\end{enumerate}
\end{corollary}

\section{Prime Polynomials\label{prime sec}}

We assume $q\neq 1$ throughout this section and apply results from section %
\ref{Normal Ore Sec} with $R=k\left[ x\right] $ and $t=y$. \ Then Theorem %
\ref{Factor Theorem} applies to quadratic forms with $b=e=0$ and $c=1$ in
equation \ref{QF Eq}. \ Remark \ref{Center Desc} and Part 3 of Lemma \ref%
{Normal Lemma} implies that an irreducible quadratic form of this type is
normal if and only if $d=0$ and $q=-1$. \ This turns out to be the most
interesting case.

The possible prime polynomials in a quantum plane are limited by \cite[%
Theorem 12]{CP}. \ If $q$ is not a root of unity, then a prime polynomial in 
$\mathcal{O}_{q}\left( k^{2}\right) $ must be a scalar multiple of $x$ or of 
$y$. \ If $q$ is a root of unity, then a prime polynomial in $\mathcal{O}%
_{q}\left( k^{2}\right) $ is either central and irreducible or a scalar
multiple of $x$ or of $y$. \ Thus, by Remark \ref{Center Desc}, there can
only be prime quadratic forms in $\mathcal{O}_{q}\left( k^{2}\right) $ if $%
q=-1$. \ Theorem \ref{Primes in QP Theorem} extends \cite[Theorem 14]{CP},
which only considered homogeneous quadratic forms with real coefficients.

\begin{theorem}
\label{Primes in QP Theorem}Let $F=ax^{2}+bxy+cy^{2}+f$ be an irreducible
quadratic form in $\mathcal{O}_{q}\left( k^{2}\right) $. \ Then $F$ is prime
if and only if $b=0$ and $q=-1$.
\end{theorem}

\begin{proof}
If $b\neq 0$ then $F\cdot x=x\left( ax^{2}+qbxy+q^{2}cy^{2}+f\right) $, $%
F\nmid x$, and $F\nmid \left( ax^{2}+qbxy+q^{2}cy^{2}+f\right) $. \ The
comments preceding the theorem explain why $q$ must be $-1$. \ We must prove
that if $q=-1$ and $F=ax^{2}+cy^{2}+f$ is irreducible in $\mathcal{O}%
_{-1}\left( k^{2}\right) $ then $F$ is prime. \ If $a=0$ or $c=0$ then this
follows from \cite[Theorem 13]{CP}. \ We easily pass to the case $c=1$. \ We
assume $F$ is irreducible but not prime in $\mathcal{O}_{-1}\left(
k^{2}\right) $ and show this leads to a contradiction.

We match our notation to Theorem \ref{Factor Theorem}. \ Set $R=k\left[ x%
\right] $ and let $Q$ be the quotient ring of $R$. \ Since $q=-1$ we know $F$
is central by Remark \ref{Center Desc} with $n=2$. \ Then $F$ is reducible
in $Q\left[ y;\sigma \right] $ and there exists $w\in Q\left[ y;\sigma %
\right] $ such that $ax^{2}+f=-w\sigma \left( w\right) $ by Theorem \ref%
{Factor Theorem}. \ Since $w\in Q$ we may find relatively prime $r,s\in k%
\left[ x\right] $ such that $w=rs^{-1}$. \ Clearing denominators gives
equation \ref{mult eq 1}. 
\begin{equation}
\left( ax^{2}+f\right) \sigma \left( s\right) s=-\sigma \left( r\right) r
\label{mult eq 1}
\end{equation}

If $ax^{2}+f$ is irreducible then equation \ref{mult eq 1} gives $\left(
ax^{2}+f\right) |r$, hence $r=\left( ax^{2}+f\right) r_{0}$ for some $%
r_{0}\in k\left[ x\right] $. \ We substitute into equation \ref{mult eq 1}
and cancel common factors to obtain equation \ref{mult eq 1a}, which implies 
$\left( ax^{2}+f\right) |s$. \ This is a contradiction since $r$ and $s$ are
relatively prime.%
\begin{equation}
-\sigma \left( s\right) s=\sigma \left( r_{0}\right) r_{0}\left(
ax^{2}+f\right)  \label{mult eq 1a}
\end{equation}

Thus $ax^{2}+f=\left( ax+\epsilon \right) \left( x+\frac{f}{\epsilon }%
\right) $ for some $\epsilon \in k$. \ Expanding and equating coefficients
gives $\epsilon ^{2}=-af$. \ Equation \ref{mult eq 1} implies $\deg r=\deg
s+1$ so by expanding and equating the leading coefficients we find $a=\nu
^{2}$ for some nonzero $\nu \in k$. \ We have proved $\epsilon ^{2}=-af$ and 
$\Delta _{-1}F=4\nu ^{2}$ which contradicts Theorem \ref{Quad Forms in QP}.
\end{proof}

More details about primes in quantum planes are provided in \cite{CP}. \ In
the rest of this section we derive analogous results for quantized Weyl
algebras.

\begin{example}
\label{variables}Suppose $q\neq -1$. \ Then $x\nmid q^{2}xy+1+q$ in $%
\mathcal{A}_{1}^{q}\left( k\right) $ by Corollary \ref{Quad Form Coro}. \
Thus $x$ is not prime in $\mathcal{A}_{1}^{q}\left( k\right) $ since $%
y^{2}x=\left( q^{2}xy+1+q\right) y$. \ Similarly, $y$ is not prime in $%
\mathcal{A}_{1}^{q}\left( k\right) $.
\end{example}

\begin{example}
\label{Normal Polys}The element $u=\left( q-1\right) xy+1$ is normal in $%
\mathcal{A}_{1}^{q}\left( k\right) $ since $ux=qxu$ and $yu=quy$. \
Moreover, $u$ is a prime polynomial since $\mathcal{A}_{1}^{q}\left(
k\right) /\left( u\right) \cong k\left[ x,x^{-1}\right] $ by \cite[%
Proposition 8.2]{Good}. \ A formula for the the coefficients of $u^{m}$ in $%
\mathcal{A}_{1}^{q}\left( k\right) $ when $q$ is a root of unity is provided
in \cite[13.5]{GL}. \ If $q$ is not a root of unity then a straightforward
induction argument shows $u^{m}=\sum_{i=0}^{m}\mu _{m,i}x^{i}y^{i}$ for each 
$m\in \mathbb{N}$ where $\mu _{m,0}=1$ and the remaining coefficients are
defined recursively by equation \ref{recursive formula}. 
\begin{equation}
\mu _{m,j+1}=\frac{q^{m}-q^{j}}{1+q+\cdots +q^{j}}\mu _{m,j}
\label{recursive formula}
\end{equation}
\end{example}

\begin{theorem}
\label{Main Theorem}Suppose $q\neq 1$ and let $u=\left( q-1\right) xy+1$.

\begin{enumerate}
\item If $q$ is not a root of unity then a polynomial in $\mathcal{A}%
_{1}^{q}\left( k\right) $ is prime if and only if it is a scalar multiple of 
$u$.

\item Suppose $P$ is an irreducible polynomial belonging to $k\left[ x\right]
$ or $k\left[ y\right] $. \ If $P$ is central in $\mathcal{A}_{1}^{q}\left(
k\right) $ then $P$ is prime in $\mathcal{A}_{1}^{q}\left( k\right) $.

\item Choose arbitrary $a,b,c,f\in k$.

\begin{enumerate}
\item If $F=bxy+f$ then $F$ is prime if and only if $F=fu$ and $f\neq 0$.

\item If $F=ax^{2}+cy^{2}+f$ and $q=-1$\ then $F$ is prime if and only if $F$
is irreducible.
\end{enumerate}
\end{enumerate}
\end{theorem}

\begin{proof}
\textbf{Part 1. }\ Suppose $P\left( x,y\right) \in \mathcal{A}_{1}^{q}\left(
F\right) $ is prime. \ Given scalars $\lambda ,\mu \in k$ we form a new
polynomial, denoted $P\left( \lambda x,\mu y\right) $, by replacing $x$ and $%
y$ by the scalar multiples $\lambda x$ and $\mu y$, respectively. \ First we
prove equations \ref{Prime 1} and \ref{Prime 2} hold for all $i,j$ such that 
$x^{i}y^{j}$ is in the support of $P\left( x,y\right) $.%
\begin{equation}
P\left( x,y\right) =q^{i-j}P\left( x/q,qy\right)   \label{Prime 1}
\end{equation}%
\begin{equation}
P\left( x,y\right) =q^{j-i}P\left( qx,y/q\right)   \label{Prime 2}
\end{equation}
A straightforward calculation gives $\left( P\left( x,y\right) \right)
\left( u\right) =\left( u\right) \left( P\left( x/q,qy\right) \right) $. \
If $P\left( x,y\right) |u$ then $P\left( x,y\right) $ is a scalar multiple
of $u$ since $u$ is irreducible.  Equations \ref{Prime 1} and \ref{Prime 2}
clearly hold for $u$. \ We are left with $P\left( x,y\right) |P\left(
x/q,qy\right) $ since $P$ is prime. \ By considering degrees, we have $%
P\left( x,y\right) =\alpha P\left( x/q,qy\right) $ for some nonzero $\alpha
\in k$. \ Equating coefficients gives $\alpha =q^{i-j}$, which proves
equation \ref{Prime 1} holds. \ A similar argument can be used for equation %
\ref{Prime 2}.

Example \ref{variables} implies $P$ is not a scalar multiple of $x$ or $y$
since $q$ is not a root of unity. \ Moreover, $P$ is irreducible so there
must be terms $x^{a}y^{0}$ and $x^{0}y^{b}$ in the support of $P$ for some $%
a,b\in \mathbb{N}$. \ Comparing coefficients in equations \ref{Prime 1} and %
\ref{Prime 2} proves $a=b=0$ and if $x^{i}y^{j}$ is in the support of $P$
then $i=j$. \ Therefore we may write $P$ as in equation \ref{p form} for
some $\lambda _{0},\ldots ,\lambda _{n}\in k$ with $\lambda _{n}\neq 0$.%
\begin{equation}
P=\sum_{i=0}^{n}\lambda _{i}x^{i}y^{i}  \label{p form}
\end{equation}%
A straightforward calculation gives $Px=xG$ where $G$ is given in equation %
\ref{g form}.%
\begin{equation}
G=\sum_{i=0}^{n}q^{i}\lambda _{i}x^{i}y^{i}+\sum_{i=1}^{n}\left( 1+q+\cdots
+q^{i-1}\right) \lambda _{i}x^{i-1}y^{i-1}  \label{g form}
\end{equation}%
Now $P$ is prime and $P\nmid x$ so $P|G$. \ Since $P$ and $G$ have the same
degree it is easy to see $P=q^{-n}G$ and $\left( 1+q+\cdots +q^{j}\right)
\lambda _{j+1}=\left( q^{n}-q^{j}\right) \lambda _{j}$ for all $j$. \ Since
the coefficients of $P$ and $u$ both satisfy the recursive formula in
equation \ref{recursive formula} we must have $P=\lambda _{0}u^{n}$ for some 
$n\in \mathbb{N}$. \ But $P$ is irreducible so $n=1$, that is, $P=\lambda
_{0}u$.

\noindent \textbf{Part 2. \ }This can be proved in the same way as \cite[%
Theorem 13]{CP}.

\noindent \textbf{Part 3. \ }In part 3(a) we have $F=bxy+f$. \ A
straightforward calculation gives $Fx=xG$ with $G=\left( qb\right) xy+\left(
b+f\right) $. \ Then $F|xG$ and $F\nmid x$ so $F|G$. \ This implies $%
G=\lambda F$ for some $\lambda \in k$ since $F$ and $G$ have the same
degree. \ Equating coefficients yields $q=\lambda $, $b=\left( q-1\right) f$%
, and $F=fu$, as desired. \ On the other hand if $f\neq 0$ and $F=fu$ then $%
F $ is prime since $u$ is prime.

In part 3(b) we must prove that if $q=-1$ and $F=ax^{2}+cy^{2}+f$ is
irreducible then $F$ is prime. \ If $a=0$ or $c=0$ then this follows from
part 2 and Remark \ref{Center Desc}. \ We easily pass to the case $c=1$. \
We assume $F$ is irreducible but not prime in $\mathcal{A}_{1}^{-1}\left(
k\right) $ and show this leads to a contradiction.

We match our notation to Theorem \ref{Factor Theorem} with $R=k\left[ x%
\right] $. \ Then, by Equation \ref{eulerian derivative}, $\delta $ is an
inner $\sigma $-derivation on $Q$ and there is an isomorphism $Q\left[
y;\sigma ,\delta \right] \cong Q\left[ t;\sigma \right] $ which is the
identity map on $Q$ and sends $y+\left( 2x\right) ^{-1}$ to $t$ by \cite[%
Lemma 1.5]{Good}. \ The image of $F=y^{2}+ax^{2}+f$ in $Q\left[ t;\sigma %
\right] $ is $t^{2}-P$ where $P=-\left( ax^{2}+\left( 2x\right)
^{-2}+f\right) $. \ By Theorem \ref{Factor Theorem} there exists $w\in Q$
such that $P=\sigma \left( w\right) w$.

Since $w\in Q$ we may find relatively prime $r,s\in k\left[ x\right] $ such
that $w=rs^{-1}$. \ Multiply both sides of $P=\sigma \left( w\right) w$ by $%
4x^{2}\sigma \left( s\right) s$ and expand to obtain equation \ref{mult eq2}%
. 
\begin{equation}
-\left( 4ax^{4}+4fx^{2}+1\right) \sigma \left( s\right) s=4x^{2}\sigma
\left( r\right) r  \label{mult eq2}
\end{equation}

Set $G=4ax^{4}+4fx^{2}+1$. Then $G$ has an irreducible factor in $k\left[ x%
\right] $ so there exist monic $P_{1},P_{2}\in k\left[ x\right] $ such that $%
P_{1}$ is irreducible and $G=4aP_{1}P_{2}$. \ We consider all possibilities
for $P_{1}$ and derive a contradiction in each case.

\begin{description}
\item[Case 1] Suppose $P_{1}=x^{2}+\alpha x+\beta $ for some nonzero $\alpha
,\beta \in k$. \ Then $P_{2}=x^{2}-\alpha x+\beta $. \ If we set $\tau =%
\frac{1}{2\beta }$ and $\omega =\alpha \tau $ then it is easy to show $\tau
^{2}=a$ and $\tau -f=\omega ^{2}$. \ By part 2(b) of Corollary \ref{Quad
Form Coro} we have $F$ is reducible in $\mathcal{A}_{1}^{-1}\left( k\right) $%
, which is a contradiction. \ \textbf{Note:} this argument shows the
coefficient of $x$ must be zero in any quadratic factor of $G$.

\item[Case 2] Suppose $P_{1}$ has degree 3. \ Then $P_{2}$ has degree 1 and
we pass to Case 3.

\item[Case 3] Suppose $P_{1}=x+\gamma $ for some nonzero $\gamma \in k$. \
Then $x-\gamma $ is also a factor of $G$ so $P_{2}=\left( x-\gamma \right)
\left( x^{2}+\lambda \right) $ for some nonzero $\lambda \in k$ by the note
we made in case 1. \ We pass to case 4 if $x^{2}+\lambda $ is irreducible. \
Otherwise $x^{2}+\lambda =\left( x+\rho \right) \left( x-\rho \right) $ for
some $\rho \in k$. \ We set $\tau =\frac{1}{2\gamma \rho }$ and $w=\tau
\left( \rho +\gamma \right) $ and argue as in case 1.

\item[Case 4] Suppose $P_{1}=x^{2}+\lambda $ for some nonzero $\lambda \in k$%
. \ Then equation \ref{mult eq2} implies $P_{1}|r$, hence $r=P_{1}r_{1}$ for
some $r_{1}\in k\left[ x\right] $. \ We substitute into equation \ref{mult
eq2} and cancel common factors to obtain equation \ref{mult eq3}. 
\begin{equation}
-P_{2}\sigma \left( s\right) s=4x^{2}\sigma \left( r_{1}\right) r_{1}P_{1}
\label{mult eq3}
\end{equation}%
Since $r$ and $s$ are relatively prime, equation \ref{mult eq3} implies $%
P_{1}|P_{2}$. \ Clearly $P_{1}=P_{2}$ in this case. \ Expand $G=4aP_{1}P_{2}$
and compare coefficients to see that $a=f^{2}$. \ By part 2(a) of Corollary %
\ref{Quad Form Coro} we have $F$ is reducible in $\mathcal{A}_{1}^{-1}\left(
k\right) $, which is a contradiction.

\item[Case 5] Suppose $P_{1}=\frac{1}{4a}G$. \ Then equation \ref{mult eq2}
implies $G|r$, hence $r=Gr_{2}$ for some $r_{2}\in k\left[ x\right] $. \ We
substitute into equation \ref{mult eq2} and cancel common factors to obtain
equation \ref{mult eq4}.%
\begin{equation}
-\sigma \left( s\right) s=4x^{2}\sigma \left( r_{2}\right) r_{2}G
\label{mult eq4}
\end{equation}%
But equation \ref{mult eq4} implies $G|s$, which is a contradiction since $r$
and $s$ are relatively prime.
\end{description}
\end{proof}

\paragraph{Acknowledgements}

The first author, Mrs. Holtz, was funded by an "Undergraduate Student and
Faculty Collaborative Research Project" provided by the University of
Wisconsin\ Oshkosh. \ Both authors would like to thank the referee for
helpful comments.


\begin{thebibliography}{99}
\bibitem{AVDVO} M. Awami, M. Van den Bergh, and F. Van Oystaeyen, Note on
Derivations of Graded Rings and Classification of Differential Polynomial
Rings, \textit{Bull. Soc. Math. Belg. S\'{e}r. A} \textbf{40 }(1988),
175-183.

\bibitem{Bueso} J. L. Bueso, J. G\'{o}mez-Torrecillas, and A. Verschoren, 
\textit{Algorithmic methods in non-commutative algebra}, in: Mathematical
Modelling: Theory and Applications, 17. Kluwer Academic Publishers,
Dordrecht, 2003.

\bibitem{Co} S. C. Coutinho, \textit{A Primer of Algebraic }$D$-\textit{%
modules}, in: London Mathematical Society Student Texts 33, Cambridge
University Press, Cambridge, 1995.

\bibitem{CP} R. Coulibaly and K. Price, Factorization in Quantum Planes,\ 
\textit{Missouri Journal of Mathematical Sciences} \textbf{18(3) }(2006).

\bibitem{Good} K. R. Goodearl, Prime Ideals in Skew Polynomial Rings and
Quantized Weyl Algebras, \textit{Journal of Algebra} \textbf{150} (1992),
324-377.

\bibitem{GL} K. R. Goodearl and E. S. Letzter, \textit{Prime Ideals in Skew
and }$q$\textit{-Skew Polynomial Rings}, in: Memoirs of the American
Mathematical Society 109\ (521), Providence, 1994.

\bibitem{GW} K. R. Goodearl and R. B. Warfield Jr., \textit{An Introduction
to Noncommutative Noetherian Rings},\ second edition, in: London
Mathematical Society Student Texts, 61. Cambridge University Press,
Cambridge, 2004.

\bibitem{Jordan} D. A. Jordan, Normal Elements of Degree One in Ore
Extensions, \textit{Communications in Algebra }\textbf{30(2)} (2002),
803-807.

\bibitem{McRob} J. C. McConnell and J. C. Robson, \textit{Noncommutative
Noetherian Rings},\ Wiley, Chichester, England, 1987.

\bibitem{Price} K. Price, A Domain Test for Lie Color Algebras, \textit{%
Journal of Algebra and its Applications }\textbf{7(1) }(2008), 81--90.
\end{thebibliography}
\end{document}